\documentclass{amsart}

\usepackage{amsmath,amssymb,amsbsy,amsfonts,latexsym,amsthm,amscd,amsxtra,
amssymb}
\usepackage[all]{xy}

\newtheorem{theorem}{Theorem}[section]
\newtheorem{lemma}[theorem]{Lemma}
\newtheorem{proposition}[theorem]{Proposition}
\newtheorem{corollary}[theorem]{Corollary}

\theoremstyle{definition}
\newtheorem{definition}[theorem]{Definition}

\theoremstyle{remark}
\newtheorem{remark}[theorem]{Remark}

\newtheorem*{Pf}{Proof}
\newenvironment{Proof}{\begin{Pf} \begin{upshape}} {\end{upshape} \qed\end{Pf}}

\numberwithin{equation}{section}

\newcommand\beqa[1]{ \begin{eqnarray} \label{#1}}
\newcommand{\eeqa}{ \end{eqnarray} }
\newcommand{\beqano}{ \begin{eqnarray*} }
\newcommand{\eeqano}{ \end{eqnarray*} }

\newcommand{\T}{ {\mathbb T}   }

\newcommand{\R}{ {\mathbb R}   }

\renewcommand \a {\alpha}
\newcommand \e {\varepsilon }
\newcommand \x {\xi}
\renewcommand \b  {\beta}

\newcommand \m {\mu}

\newcommand \om {\omega}

\newcommand \g {\gamma}

\newcommand \G {\Gamma}

\renewcommand \l {\lambda}
\renewcommand \L {\Lambda}

\newcommand \cA {{\mathcal A}}
\newcommand \cC {{\mathcal C}}

\newcommand \cF {{\mathcal F}}
\newcommand \cL {{\mathcal L}}
\newcommand \cM {{\mathcal M}}
\newcommand \cN {{\mathcal N}}

\newcommand \cS {{\mathcal S}}

\newcommand \fg {{\mathfrak g}}
\def\ie{\hbox{\it i.e.\ }}

\newcommand \dpr {\partial}

\newcommand \rH {{\rm H}}
\newcommand \rT {{\rm T}}
\newcommand \supp {{\rm supp\ }}
\newcommand \rank {{\rm rank\ }}

\begin{document}

\title{On the integrability of Tonelli Hamiltonians}

\author{Alfonso Sorrentino}
\address{CEREMADE, UMR du CNRS 7534\\ Universit\'e Paris-Dauphine\\
Place du Mar\'echal De Lattre De Tassigny\\ 75775 Paris Cedex 16\\ France.}
\curraddr{Department of Mathematics and Mathematical Statistics\\
University of Cambridge\\ Wilberforce Road\\ Cambridge CB3 0WB\\ United Kingdom.}
\email{A.Sorrentino@dpmms.cam.ac.uk}
\urladdr{http://www.dpmms.cam.ac.uk/~as998}
\thanks{I am grateful to Patrick Bernard for his interest in this work and for his invaluable help. I would also like to thank the anonimous referee for several helpful comments and suggestions, which contributed to simplify some proofs.}

\subjclass[2000]{Primary 37J50, 37J35, 37J15; Secondary 53D12, 53D25}


\begin{abstract}
In this article we  discuss a weaker version of Liouville's theorem on the integrability of Hamiltonian systems.  
We show that in the case of Tonelli Hamiltonians
the {involution} hypothesis on the integrals of motion can be completely dropped and still interesting information on the dynamics of the system can be deduced. Moreover, we prove that on the $n$-dimensional torus this weaker condition implies classical integrability in the sense of Liouville.
The main idea of the proof consists  in relating the existence of independent integrals of motion of a Tonelli Hamiltonian to the ``size'' of its
Mather and Aubry sets. As a byproduct we point out
the existence of ``non-trivial'' common invariant sets  
for all Hamiltonians that Poisson-commute with a Tonelli one. \end{abstract}

\maketitle


\section{Introduction}\label{intro}

A classical result in the study of Hamiltonian systems is what is generally called {\it Liouville's Theorem} (or  {\it Arnol'd-Liouville's Theorem}, see \cite{Arnoldbook}), which is concerned with the {\it integrability} of these systems, \ie the possibility of expressing their solutions in a closed form. Actually, in the case of Hamiltonian systems this notion assumes a more precise connotation, sometimes specified as {\it integrability in the sense of Liouville}: it refers to the existence of  a regular foliation of the phase space by invariant Lagrangian submanifolds, which are diffeomorphic to tori and on which the dynamics is conjugate to a rigid rotation.
Liouville's theorem provides sufficient conditions for the existence of such a foliation. 
Let us state  Liouville's result in a more precise way.
Consider a Hamiltonian system with $n$ degrees of freedom, given by a Hamiltonian $H: V \rightarrow \R$ defined on a $2n$-dimensional symplectic manifold $(V,\omega)$ and denote by $\{\cdot,\cdot\}$ the associated {\it Poisson bracket}, defined as follows: if $f,g \in C^1(M)$, then $\{f,g\}=\om(X_f,X_g)=df\cdot X_g$, where $X_f$ and $X_g$ denote the Hamiltonian vector fields associated to $f$ and $g$ (see for instance \cite{Arnoldbook}).  
An important role in the study of the dynamics  is played by the functions $F: V \rightarrow \R$ that are in  {\it involution} with the Hamiltonian, \ie whose Poisson bracket $\{H,F\}\equiv 0$ on $V$ (equivalently we can say that $H$ and $F$ {\it Poisson-commute}). Such functions, whenever they exist, are called {\it integrals of motion} (or {\it first integrals}) of $H$. It is quite easy to check that the condition of being in involution  is equivalent to ask that $F$ is constant along the orbits of the Hamiltonian flow of $H$ and {\it vice versa}; moreover, this implies that the associated Hamiltonian vector fields $X_H$ and $X_F$ commute (see Lemma \ref{lemma1}).
Liouville's theorem relates the integrability of a given Hamiltonian system to the existence of ``enough''  integrals of motion.\\
 
\noindent {\bf Theorem [Liouville]} {\it Let $(V,\omega)$ be a symplectic manifold with $\dim V=2n$ and let $H: V \longrightarrow \R$ be a Hamiltonian. Suppose that there exist $n$ integrals of motion $F_1, \ldots, F_n : V \longrightarrow   \R$ such that:
\begin{itemize}
\item[{\rm i)}] $F_1,\ldots, F_n$ are independent at each point of $V$, \ie  $dF_1,\ldots, dF_n$ are linearly independent as vectors at each point of $V$;
\item[{\rm ii)}] $F_1,\ldots, F_n$ are pairwise in involution, \ie $\{F_i,F_j\}=0$ for all $i,j=1,\ldots n$.
\end{itemize}
Then, all non-empty level sets $\Lambda_{a} := \{F_1=a_1, \ldots, F_n=a_n\}$, where $a=(a_1,\ldots,a_n)$ is a vector in $\R^n$, are invariant under the Hamiltonian flow of $H$ and each of their connected components is a smooth Lagrangian submanifold. Moreover, such Lagrangian submanifolds are diffeomorphic to $n$-dimensional tori and the dynamics on them are conjugate to rigid rotations on $\T^n$.\\

}

\noindent See for instance \cite[Section 49]{Arnoldbook} for a proof of this theorem.

\begin{remark}\label{Remintro}
It turns out from the above theorem  
that having $n$ independent integrals of motion, pairwise in involution, is a very strong assumption with  significant implications on the dynamics of the system and the topology of the underlying configuration space.
Observe, in fact,  that while the invariance of these $\L_a$'s 
simply follows from $F_i$ being integrals of motion (and Lemma \ref{lemma1}), the fact that these submanifolds are Lagrangian, that they are diffeomorphic to tori and that the motion on each of them is conjugate to a rigid rotation, strongly relies on  these integrals of motion being  pairwise in involution. 
\end{remark}

Is it possible to weaken the assumptions in Liouville's theorem? It is easy to see that condition i) can be replaced by assuming that the integrals of motion are not independent everywhere, but in an open dense set of the space. Then,  the same conclusions will continue to hold for all level sets that are not in the {\it singular set} of the integrals of motion.   
On the other hand,
in the light of Remark \ref{Remintro}, condition ii) is definitely less easy to weaken or to drop. In \cite[Theorem 4.1]{Fomenko}, for instance, Fomenko and Mishchenko proved a {\it non-commutative} version of Liouville's theorem, in which they consider the case in which  only a subfamily of the integrals of motion are pairwise in involution; as a counterbalance, they need to require that extra independent integrals of motion exist. As we shall point out in Remark \ref{remnoexist}, this cannot happen for a Tonelli Hamiltonian (unless one excludes the Mather and Aubry sets from the domain of independence of the integrals).
Observe moreover that the statement in \cite{Fomenko} would become trivial if one assumed that none of the integrals of motion were in involution; in fact, this would require the existence of $2n$ integrals of motion, exactly as many as the dimension of the phase space.\\

In this work we would like to address the following question: what happens when the involution hypothesis on the integrals of motion is completely dropped. Let us introduce the following   definition.

\begin{definition}[{\bf Weak integrability}]\label{weaklyint}
Let $H$ be a $C^2$ Hamiltonian on $(V,\omega)$ and let $\dim V =2n$.  $H$ is called {\it weakly integrable} if there exist $n$ integrals of motion, which are $C^2$  and independent at each point of $V$.
\end{definition}

\begin{remark}
({\it i}) As far as the regularity of the integrals of motion is concerned, we need them to be  $C^2$ (see for instance the proof of Lemma \ref{weakLiouville}). However, with some technical modifications, our proof could actually work assuming that the integrals of motion are only $C^{1,1}$.
({\it ii}) As already pointed out, we are assuming that these integrals of motion are independent everywhere. One could actually require them to be independent almost everywhere, but this would require extra assumptions on their action-minimizing sets. See for instance Remark \ref{remnoexist}.
\end{remark}

\noindent 
Certainly, Liouville integrable Hamiltonians are weakly integrable. As we shall prove in Appendix \ref{App1} (see Proposition \ref{propgruppilie}), other interesting examples of weakly integrable Hamiltonians, which are not necessarily integrable in the classical sense, are provided by geodesic flows associated to left-invariant Riemannian metrics on compact Lie groups.\\

In the following we shall show that in the case of {\it Tonelli Hamiltonians} - that is Hamiltonians defined on the cotangent space of a compact connected finite-dimensional manifold, which are strictly convex and superlinear in each fiber (see section \ref{secmain}) - one can deduce interesting information on the dynamics of weakly integrable Hamiltonian systems. 
More precisely (see section \ref{proofweakLiouville} for a more precise statement):\\

\noindent {\bf Theorem \ref{teolio} (Weak Liouville's theorem).}
{\it If  $H: \rT^*M \rightarrow \R$ is a weakly integrable Tonelli Hamiltonian, then for each cohomology class $c\in \rH^1(M;\R)$ there exists a unique smooth invariant Lagrangian graph $\L_c$ with cohomology class $c$, on which the dynamics is recurrent, \ie each orbit returns infinitely many often close to its initial data}.\\  

\begin{remark}
({\it i}) The smoothness of these Lagrangian graphs is related to the smoothness of the integrals of motions $F_i$ (see Lemma \ref{lemma123} and Remark \ref{remarkregularity}). In particular, if all  integrals of motion are $C^k$, with $k\geq 2$, then these graphs will be $C^{k}$ too. 
({\it ii}) Observe that this also implies that the associated Hamilton-Jacobi equation $H(x,c+du)=\mbox{\it const.}$ admits a (unique) smooth solution for each cohomology class $c \in \rH^1(M;R)$.
\end{remark}

Moreover, we shall prove that in some cases (essentially on the torus) this notion of integrability coincides with the classical one.\\

\noindent 
\noindent {\bf Theorem \ref{teolio2}.}
{\it 
Assume $\dim \rH^1(M;\R)=\dim M$ and 
let $H: \rT^*M \rightarrow \R$ be a weakly integrable Tonelli Hamiltonian with integrals of motion $F_1,\ldots, F_n$. Then the system is integrable in the sense of Liouville. In particular, $M$ is diffeomorphic to $\T^n$.}\\

\noindent Observe that this implies that the integrals of motion are pairwise in involution everywhere, although we had not assumed it {\it a-priori}.\\

The main idea behind our approach consists in studying how the existence of independent integrals of motion of a Tonelli Hamiltonian $H$ relates to the structure of its {\it action-minimizing sets}, namely its {\it Mather } and {\it Aubry sets} (see section \ref{secmain} for a definition); moreover, using the symplectic structure of the Aubry set 
(see (\ref{Aubry}) and Remark \ref{Remark11}) we shall be able to recover the involution hypothesis at least on these sets. The key properties that we are going to use can be summarized as follows:
\begin{itemize}
\item {the Mather and Aubry sets are invariant under the flow of any integral of motion of $H$} [Lemma \ref{maintheorem}];
\item {the existence of $k$ independent integrals of motion implies that the ``size'' of each Mather and Aubry set is bigger or equal than $k$} [Proposition \ref{Prop100}]; 
\item {the integrals of motion are locally in involution on the Mather and Aubry sets}  [Proposition \ref{step2}].
 \end{itemize}

As a byproduct of this discussion, we shall also point out some results about the existence of ``non-trivial'' common invariant  sets (or measures) for Poisson-commuting Hamiltonians. It will turn out, in fact, that  \noindent{\it  if $H$ is a Tonelli Hamiltonian, then there exists a family of  ``non-trivial'' compact invariant subsets, which are also invariant under the flows of ALL Hamiltonians that Poisson-commute with $H$}  [Proposition \ref{Prop1}].\\  

The paper is organized as follows. In Section \ref{secmain}, after having  briefly recalled what Mather's theory is about and introduced the notion of Mather and Aubry sets, we shall discuss the relation between these sets and the integrals of motion of the Hamiltonian. In Section \ref{proofweakLiouville} we shall prove the main results announced in this Introduction. In Section \ref{secproof} we shall prove a fundamental lemma, stated in Section \ref{secmain}, about the symplectic invariance of the Mather and Aubry sets. To conclude, in Appendix \ref{App1} we shall discuss some examples of weakly integrable Tonelli Hamiltonian systems, namely the geodesic flows associated to left-invariant Riemannian metrics on compact Lie groups.


\section{Aubry-Mather sets and Integrals of motion}\label{secmain}

In this section we would like to discuss the relation between the Aubry-Mather sets of a Tonelli Hamiltonian and its integrals of motion, with particular attention to how their flows act on these sets.
Before entering into details,  let us try to describe ``in a nutshell'' what {\it Mather's theory} is about, in order to help an unaware reader  to understand the concepts and the results that we shall be dealing with, although we refer her or him to \cite{Fathibook,Mather91, Mather93} for more exhaustive presentations of the material.

Mather's theory  
consists in a variational approach to the study of convex Lagrangian systems, called {\it Tonelli Lagrangians}, with particular  attention to their {\it action-minimizing invariant probability measures}  and their {\it action-minimizing orbits}. 
Let $M$ be a compact and connected smooth $n$-dimensional manifold without boundary. Denote by ${\rm T}M$ its tangent bundle and ${\rm T}^*M$ the cotangent one and denote points in ${\rT M}$ and ${\rT^*M}$
respectively by $(x,v)$ and $(x,p)$. We shall also assume that the cotangent bundle $\rT^*M$ is equipped with the canonical symplectic structure, which we shall denote $\omega$.
A {\it Tonelli Lagrangian} is a $C^2$ function $L: \rT M \rightarrow \R$, which is strictly convex and uniformly superlinear in the fibers; in particular this Lagrangian defines a flow on ${\rT M}$, known as {\it Euler-Lagrange flow}. 
Instead of considering just this Lagrangian $L$, John Mather \cite{Mather91} proposed to consider a family of modified Tonelli Lagrangians given by 
$L_{\eta}(x,v)=L(x,v)-\langle \eta(x),v\rangle$, where $\eta$ is a closed $1$-form on $M$. These Lagrangians, in fact, have all the same Euler-Lagrange flow as $L$, but different action-minimizing orbits and measures, according to the cohomology class of $\eta$.
In this way for each $c \in{\rm H}^1(M;\R)$, if we choose $\eta$ to be any smooth closed $1$-form on $M$ with cohomology class $[\eta]=c$ and we consider the Lagrangian $L_{\eta}$, it is possible to 
define two compact invariant subsets of ${\rm T}M$:
\begin{itemize}
\item $\widetilde{\cM}_c(L)$, the {\it Mather set of cohomology class $c$}, given by the union of the supports of all invariant probability measures that minimize the action of $L_{\eta}$ ($c$-{\it action minizimizing measure} or {\it Mather's measures of cohomology class} $c$); 
\item $\widetilde{\cA}_c(L)$, the {\it Aubry set of cohomology class $c$}, given by the union of  all {\it regular global minimizing curves}  of $L_{\eta}$ (or $c$-{\it regular minimizers}); see \cite{Mather93, Fathibook} for a precise definition. 
\end{itemize} 

\noindent These sets are such that $\widetilde{\cM}_c(L) \subseteq \widetilde{\cA}_c(L)$ and moreover one of their most important features is that they are graphs over $M$ ({\it Mather's graph theorem} \cite{Mather91, Mather93}); namely, if $\pi: {\rm T}M \rightarrow M$ denotes the canonical projection along the fibers, then $\pi|{\widetilde{\cA}_c(L)}$ is injective and its inverse $\big(\pi|\widetilde{\cA}_c(L)\big)^{-1}\!\!\!: \pi\big(\widetilde{\cA}_c(L)\big) \longrightarrow \widetilde{\cA}_c(L)$ is Lipschitz. \\
In the following we would like to consider the Hamiltonian setting rather than the Lagrangian one. 
It is a classical result that, under the above Tonelli conditions, it is possible to associate to the Lagrangian system a {Hamiltonian system} 
$H: \rT^* M \rightarrow \R$, where $H(x,p)=\sup_{v\in \rT_xM}(\langle p,v \rangle - L(x,v))$. It is easy to check that $H$ is also $C^2$, strictly convex and uniformly superlinear in each fiber: $H$ will be called a {\it Tonelli} (or sometimes {\it optical}) {\it Hamiltonian}. The important observation is that the flow on $\rT^* M$ associated to this Hamiltonian, known as the {\it Hamiltonian flow} of $H$, is conjugate  to the Euler-Lagrange flow of $L$, via the so-called {\it Legendre transform}, \ie the diffeomorphism $\cL_L: \rT M \longrightarrow \rT^*M$ defined by
$\cL_L(x,v)=(x, \frac{\dpr L}{\dpr v}(x,v))$.

Therefore, one can define the analogue of the Mather and Aubry sets in the cotanget space, simply considering $\cM^*_c(H)=\cL_L\big(\widetilde{\cM}_c(L)\big)$ and
$\cA^*_c(H)=\cL_L\big(\widetilde{\cA}_c(L)\big)$. These sets still satisfy all the properties mentioned above, including the graph theorem.
Moreover, as it was proved by Carneiro \cite{Carneiro},  these sets are contained in the energy level $\{H(x,p)=\a_H(c)\}$, where $\a_H: \rH^1(M;\R)\longrightarrow \R$ is called {\it Mather's $\a$-function} (or {\it Ma\~n\'e's critical value} or {\it effective Hamiltonian}) and $-\a_H(c)$ represents the average action of $c$-action minimizing measures; see \cite{Mather91, Fathibook} for a more precise definition. \\
Before concluding this preamble, let us recall that using Fathi's weak KAM theory \cite{Fathibook} it is possible to obtain a nice characterization of the Aubry set in terms of {\it critical subsolutions} of Hamilton-Jacobi equation. As above, let $\eta$ be a closed $1$-form with cohomology class $c$; we shall say that $u \in C^{1,1}(M)$ is an $\eta$-critical subsolution if it satisfies $H(x,\eta + d_xu)\leq \a_H(c)$ for all $x\in M$. The existence of such functions has been showed by Bernard \cite{bernardc11}.
If one denotes  by $\cS_{\eta}$ the set of $\eta$ critical subsolutions, then: 
\beqa{Aubry} {\cA}_c^*(H)=\bigcap_{u\in \cS_{\eta}} \left\{(x,\eta_x + d_xu):\; x\in M\right\}.\eeqa
\begin{remark}\label{Remark11}
Recall that in $\rT^*M$, with the standard symplectic form, there is a $1$-$1$ correspondence between Lagrangian graphs and closed $1$-forms (see for instance \cite{AnnaCannas}). Therefore, we could interpret the graphs of these critical subsolutions as {\it Lipschitz} Lagrangian graphs in $\rT^*M$. We shall see that the fact that the Aubry set can be seen as the intersection of these {special}  Lagrangian graphs (see Definition \ref{Def2} and Remark \ref{Remark1}), will play a key role in our proof.\\
\end{remark}

\vspace{5 pt}

Now that we have recalled most of the definitions that we shall need, we can start our discussion. 
Let $H$ be a Tonelli Hamiltonian on $\rT^*M$ and $F$ an integral of motion of $H$.  
As we have already mentioned, the main idea behind our approach consists in undestanding how the flow of $F$ acts on the Mather and Aubry sets of $H$. 

\begin{lemma} \label{maintheorem}
Let $H$ be a Tonelli Hamiltonian on $\rT^*M$ and $F$ an integral of motion of $H$.  
Let us denote by $\Phi_H$ and $\Phi_F$ the respective flows. Then, the following holds:
\begin{itemize}
\item[{\rm (i)}] If $\m$ is a $c$-action minimizing measure of $H$, then ${\Phi^t_F}_*\m$ is still a $c$-action minimizing measure of $H$, for each $t\in \R$,
where the lower $*$ denotes the push-forward of the measure. 
\item[{\rm (ii)}] The Mather set $\cM^*_c(H)$ and the Aubry set $\cA^*_c(H)$ are invariant under the action of $\Phi^t_F$, for each $t\in \R$ and for each $c\in \rH^1(M;\R)$. In particular, for each $t\in \R$, $\Phi^t_F$  maps each connected component of  $\cM^*_c(H)$ and $\cA^*_c(H)$ into itself.
 \end{itemize}
\end{lemma}

\vspace{10 pt}

\noindent We postpone the proof of this result to section \ref{secproof}.

\begin{remark}\label{giadim}
It is worthwhile to point out that this result can be also deduced from a  
result by Patrick Bernard  \cite[Theorem in $\S \;1.10$, page $6$]{Bernard} on the symplectic invariance of the Mather and Aubry sets.  
In fact for any fixed time $t$ the Hamiltonian flow $\Phi_F^t$ is an exact symplectomorphism that preserves $H$. Essentially, in section \ref{secproof} we shall provide a different proof of this result, but in the autonomous case.\\  
Another related result is contained in \cite{Maderna}, where the author considers the action of symmetries of the Hamiltonian, \ie $C^1$-diffeomorphisms of $M$ that preserve $H$. One can deduce from the results therein that the Mather and Aubry sets of $H$ are invariant under the action of the identity connected component of the group of such diffeomorphisms. From our point of view, these diffeomorphisms correspond to integrals of motion depending only on the $x$-variables.
\end{remark}

As recalled in Section \ref{intro}, Liouville's theorem is concerned with {\it independent} integrals of motion, \ie integrals of motion whose differentials are linearly independent, as vectors, at each point.  Let us see how the existence of independent integrals of motion relates to  the ``size'' of the Mather and Aubry sets of $H$. 
In order to make clear what we mean by ``size'' of these sets, let us introduce some notion of tangent space.
We shall call {\it generalized tangent space} to $\cM^*_c(H)$ (resp.~$\cA^*_c(H)$) at a point $(x,p)$, the set of all vectors that are tangent to curves in $\cM^*_c(H)$ (resp.~$\cA^*_c(H)$) at $(x,p)$. We shall denote it by $\rT_{(x,p)}\cM^*_c(H)$ (resp.~$\rT_{(x,p)}\cA^*_c(H))$ and we shall define its {\it rank} to be the largest number of linearly independent vectors that it contains.
Observe that 
$\rT_{(x,p)}\cM^*_c(H) \subseteq \rT_{(x,p)}\cA^*_c(H)$. In particular, if the Mather set does not contain any fixed point
(\ie $dH(x,p)\neq 0$ for all $(x,p)\in \cM^*_c(H)$), then $\rank \rT_{(x,p)}\cA^*_c(H) \geq \rank \rT_{(x,p)}\cM^*_c(H) \geq 1$; 
in fact, since these sets are invariant, the Hamiltonian vector field $X_H(x,p) \neq 0$ is tangent to them. 

\begin{proposition}\label{Prop100}
Let $H$ be a Tonelli Hamiltonian on $\rT^*M$ and suppose that there exist $k$ independent integrals of motion.
Then, $\rank \rT_{(x,p)}\cM^*_c(H) \geq k$ at all points $(x,p)\in \cM^*_c(H)$ for each $c\in \rH^1(M;\R)$.
\end{proposition}

\begin{Proof}
It follows from the fact that $\cM^*_c(H)$ is invariant under the flows of the $k$ independent integrals of motion (Lemma \ref{maintheorem}).
The linear independence of the corresponding vector fields (which are therefore tangent to this set) follows from the independence of the integrals of motion  and the non-degeneracy of the symplectic form $\omega$.
\end{Proof}

\vspace{10 pt}

Conversely, some information on the existence of independent integrals of motion of $H$ can be obtained from the
``structure'' of the Mather sets. Let us define 
$$\l(H):= \min_{c\in\rH^1(M;\R)} \min_{(x,p)\in\cM^*_c(H)} \rank T_{(x,p)} \cM^*_c(H)\,.$$ 
This quantity is clearly well-defined and, if the Mather sets do not contain any fixed point, it is bigger or equal than $1$.

\begin{corollary}\label{Prop101}
Let $H$ be a Tonelli Hamiltonian on $\rT^*M$. Then, there may exist at most $\l(H)$ independent integrals of motion of $H$.
In particular, if some  $\cM^*_{c}(H)$ contains an isolated periodic orbit, then $\l(H)=1$ and therefore all integrals of motion of $H$ are linearly dependent on $H$ on this orbit.
\end{corollary}

\begin{remark}\label{remnoexist}
({\it i}) Observe that the above results remain true if we assume that the integrals of motion are defined only locally, \ie in an open region of the phase space. In this case, we can still apply the same ideas and get information on the Mather and Aubry sets contained in this open region.
({\it ii}) From Corollary \ref{Prop101}, it results clearly that there may exist at most $n$ integrals of motion ($n=\dim M$) that are independent in a region containing  some part of Mather or Aubry sets. 
\end{remark}

Another important peculiarity of the Mather and Aubry sets, with respect to their interplay with the integrals of motion,  is that they are not only invariant under the action of their flows, but they also force the integrals of motion to Poisson-commute. In fact, using the characterization of the Aubry set in terms of critical subsolutions of Hamilton-Jacobi (see (\ref{Aubry})) and its symplectic interpretation (see Remark \ref{Remark11}), one can recover the involution property of the integrals of motion, at least locally.

\begin{proposition}\label{step2}
Let $H$ be a Tonelli Hamiltonian on $\rT^*M$ and let $F_1$ and $F_2$ be two $C^2$ integrals of motion. Then for each $c\in \rH^1(M;\R)$ we have that 
$$\{F_1,F_2\}(x,\hat{\pi}_c^{-1}\!(x))=0 \qquad \forall\, x\in \overline{{\rm Int}\big(\cA_c(H)\big)},$$ where $\hat{\pi}_c=\pi|\cA^*_c(H)$   and $\cA_c(H)=\pi\big(\cA^*_c(H)\big)$.
\end{proposition}

\begin{Proof}
From weak KAM theory we can deduce that $\cA_c^*(H) $ is contained in a Lipschitz Lagrangian graph $\L$ (see (\ref{Aubry}) and Remark \ref{Remark11}),
which is the graph of a Lipschitz closed $1$-form $\eta: M \rightarrow \rT^*M$. Observe that $\eta=\hat{\pi}_c^{-1}$ on $\cA_c(H)$. Moreover $\eta$ is differentiable almost everywhere on $M$ and at all differentiability points $x$ the tangent space $\rT_{\eta(x)} \L$ is a Lagrangian subspace. 
Let $x\in {\rm Int}\big(\cA_c(H)\big)$ be a differentiability point of $\eta$ (observe that $\eta$ is differentiable at almost every point in the interior of $\cA_c(H)$). 
Since $\cA_c^*(H)$ is invariant under the action of both flows $\Phi_{F_1}$
and $\Phi_{F_2}$, then the associated vector fields $X_{F_1}$ and $X_{F_2}$ are tangent to $\L$ on $\cA^*_c(H)$. 
Using the definition of the Poisson-bracket and the fact that $\rT_{\eta(x)} \L$ is Lagrangian, we get:
$\{F_1, F_2\}(x,\eta(x)) = \omega( X_{F_1}(x,\eta(x)), X_{F_2}(x,\eta(x))) = 0$.
By continuity, this extends to  $\overline{{\rm Int}\big(\cA_c(H)\big)}$.
\end{Proof}

\vspace{10 pt}

Before concluding this section, let us point out an interesting consequence of Lemma \ref{maintheorem}.
A natural question that someone might ask is the following.  
Suppose that we have two Poisson-commuting Hamiltonians $H$ and $F$ on ${\rT^*M}$, \ie $\{H,F\}=0$. 
 Then, {is it possible to find sets or measures that are invariant under the action of both flows?}
For instance, it is easy to check that all energy levels of $H$ are also invariant under the flow of $F$ (and {\it vice versa}). However, these ``trivial'' sets do not seem to provide a satisfactory answer; rather,
it would be interesting to show the existence of ``non-trivial'' common invariant sets, such as, for example, invariant sets that have positive codimension in a given energy level.\\ 
Lemma \ref{maintheorem} provides a positive answer to such a question in the case of {Tonelli} Hamiltonians. Actually, something more is true.

\begin{proposition}\label{Prop1}
Let $H$ be a Tonelli Hamiltonian on $\rT^*M$. Then, there exists a family of compact invariant sets of $H$, parametrized over $\rH^1(M;\R)$,  
which are invariant under the flows of ALL integrals of motion of $H$.
These sets are supported on Lipschitz Lagrangian graphs over $M$ and, if $\rH^1(M;\R)$ is not trivial, each 
energy level above Ma\~n\'e's strict critical value {\rm(}\ie the minimum of Mather's $\a$-function{)} contains at least one of them.  
\end{proposition}

\begin{Proof}
This family is given by $\cM^*_c(H)$ and $\cA^*_c(H)$ for all $c\in\rH^1(M;\R)$ (because of Lemma \ref{maintheorem}).
If $\rH^1(M;\R)$ is  not trivial, then for each energy value $E\in\R$ above the minimum of $\a_H$ (which is a convex and superlinear function), there exists at least one $c\in\rH^1(M;\R)$ such that $\a_H(c)=E$. This energy value will contain $\cM^*_c(H)$ and $\cA^*_c(H)$. 
\end{Proof}

\vspace{3 pt}

We can also deduce the following consequence. 

\begin{corollary}\label{Cor1} 
Let $H$ be a Tonelli Hamiltonian on $\rT^*M$ and $F$ an integral of motion of  $H$.
If $H$ has an invariant Lipschitz Lagrangian graph $\Lambda$ supporting an invariant measure $\m$ of full support (\ie $\supp \m = \Lambda$), then $\Lambda$ is also invariant under the flow of  $F$.
\end{corollary}

\begin{Proof} The proof follows from the fact that if such $\Lambda$ exists, then $\Lambda=\cM^*_{c_{\Lambda}}(H)=\cA^*_{c_{\Lambda}}(H)$, where ${c_{\Lambda}}$ denotes the cohomology class of $\Lambda$ (see \cite{FGS}).
\end{Proof}

\vspace{3 pt}

One could also ask about the existence of common invariant ergodic probability measures. In general $(i)$ in Lemma \ref{maintheorem} does not imply ${\Phi^t_F}_*\m=\m$ (take for instance the case of an invariant torus foliated by periodic orbits and choose $\m$ to be a measure supported on one of these orbits). However in some cases it is possible to deduce it.

\begin{corollary}\label{cor2}
Let $H$ be a Tonelli Hamiltonian on $\rT^*M$ and $F$ an integral of motion of $H$.
If $\cM^*_c(H)$ is uniquely ergodic, \ie  $H$ has a unique $c$-action minimizing measure $\m$, then  
${\Phi^t_F}_*\m= \m$ for all $t\in \R$. In other words, $\m$ is also invariant for $F$.
\end{corollary}

\begin{Proof}
From (i) in Lemma \ref{maintheorem}, it follows that ${\Phi^t_F}_*\m$ is still a $c$-action minimizing measure of $H$ for all $t\in\R$. The unique ergodicity of $\cM^*_c(H)$ implies that  necessarily ${\Phi^t_F}_*\m = \m$.
\end{Proof}

\begin{remark}
In \cite{Mane} Ricardo Ma\~n\'e  showed  that for
any given Tonelli Hamiltonian $H$, there exist residual subsets $\cS(H) \subseteq C^2(M)$ and $\cC(H)  \subseteq \rH^1(M;\R)$ such that for each $V \in \cS(H)$ and $c\in \cC(H)$, the Mather set $\cM^*_c(H+V)$ is uniquely ergodic. 
With this on mind one can deduce that for a {\it generic} Tonelli Hamiltonian $H$, all integrals of motion of  $H$ have uncountably many invariant ergodic measures in common (assuming that $\rH^1(M;\R)$ is not trivial). 
In particular, these measures are $c$-action minimizing measures for $H+V$, for some $c\in\cC(H)$ (observe that this set does not depend on $V$). Furthermore, for each $E\in\R$, infinitely many of these measures will have energy bigger than $E$ (\ie their supports will be contained in an energy level of $H+V$ bigger that $E$).
However, it is important to point out that most likely all these integrals of motion will be functionally dependent on $H+V$ at some point, since the existence of independent integrals of motion is a highly non-generic situation. 
\end{remark}

\vspace{10 pt}


\section{Weak Liouville's Theorem}\label{proofweakLiouville}

In this section we shall prove our weak version of Liouville's theorem for Tonelli Hamiltonians and the other results announced in section \ref{intro}. 
Let us start by stating our results.

\begin{theorem}[{\bf Weak Liouville's theorem}]\label{teolio}
Let $M$ be a compact connected $n$-dimensional manifold without boundary and let $H: \rT^*M \rightarrow \R$ be a weakly integrable Tonelli Hamiltonian. Then:
\begin{itemize}
\item[{\rm i)}] for each cohomology class $c\in \rH^1(M;\R)$ there exists a smooth invariant Lagrangian graph $\L_c$ with cohomology class $c$, which supports an invariant measure of full support, \ie  $\L_c=\cM^*_c(H)=\cA^*_c(H)$. In particular, the motion on each $\L_c$ is recurrent.
\item[{\rm ii)}] These Lagrangian graphs $\{\L_c\}_{c\in \rH^1(M;\R)}$ are disjoint and any other invariant Lipschitz Lagrangian graph with cohomology class $c$ must coincide with the corresponding $\L_c$ in this family.
\item[{\rm iii)}] Let us denote $\L^* := \cup_{c\in \rH^1(M;\R)} \L_c $. 
$\L^*$ is closed and 
if some $\L_c \subseteq {\rm Int}(\L^*)$, then $\L_c$ is diffeomorphic to an $n$-dimensional torus and the motion on it is conjugate to a rotation on $\T^n$; in particular, $M$ is diffeomorphic to $\T^n$. If $\L^*$ is open, then the system is integrable in the classical (Liouville) sense.
\end{itemize}
\end{theorem}

\begin{remark}\label{Remdopoteolio} 1) Observe that point i) of  Theorem \ref{teolio}  implies that the associated Hamilton-Jacobi equation $H(x,c+du)=\a_H(c)$ admit a unique smooth solution for each cohomology class $c \in \rH^1(M;R)$ ($\a_H(c)$ is the unique energy value for which a solution may exist, see for instance \cite{Fathibook}).\\
2) One could assume that the integrals of motion are not independent everywhere, but only in an open dense set of the phase space. In this case, in the light of Remark \ref{remnoexist}, the conclusions i) and ii) of Theorem \ref{teolio} will continue to hold, but only for those cohomology classes whose Aubry-Mather sets lie in this region. \\
3) As far as the uniqueness result of point ii) is concerned, it is easy to see  that these graphs are also unique in the class of smooth Lagrangian submanifolds isotopic to the zero section. It suffices to consider their graph selectors (see \cite{Chap, Siburg}) and use the same proof as in ii).\\
4) Using the disjointness of these Lagrangian graphs, one can to conclude that, as for Liouville integrable systems, $\a_H$ is strictly convex and its convex conjugate $\b_H$  is $C^1$ (Corollary \ref{Coralphabeta}). Recall that $\b_H(h) := \sup_{c\in\rH^1(M;\R)} (\langle c,h \rangle - \a_H(c))$; this function, also called Mather's $\b$-function, represents the minimal action of invariant probability measures with rotation vector $h$ (see \cite{Mather91} for a precise definition). 
\end{remark}

Now it would be interesting to understand whether or not there are cases in which this weaker notion of integrability is equivalent to the classical  one (in the sense of Liouville). As remarked in point iii) of Theorem \ref{teolio}, the union of  these Lagrangian graphs is not necessarily a foliation of the whole space (and when this happens, the manifold has to be diffeomorphic to a torus). In fact, if the dimension of $H^1(M;\R)$ is less than the dimension of $M$, this family of graphs is not sufficient to foliate $\rT^*M$ or even to have non-empty interior (for instance, think about the case in which $\rH^1(M;\R)$ is trivial). What we shall prove is that when this obstacle is removed, then the two notions coincide.

\begin{theorem}\label{teolio2}
Let $M$ be a compact connected $n$-dimensional manifold, such that $\dim \rH^1(M;\R)=\dim M$.
If $H: \rT^*M \rightarrow \R$ is a weakly integrable Tonelli Hamiltonian, then the system is integrable in the sense of Liouville. In particular, $M$ is diffeomorphic to $\T^n$.
\end{theorem}

\vspace{5 pt}

We shall split the proof of these theorems into several lemmata.
Let us start by observing that Proposition \ref{Prop100} allows us to deduce more information on the structure of the Mather and Aubry sets associated to a weakly integrable Tonelli Hamiltonian.

\begin{lemma}\label{weakLiouville}
Let $H$ be a weakly integrable Tonelli Hamiltonian on $\rT^*M$. 
Then, for each $c\in \rH^1(M;\R)$ we have that $\cM^*_c(H)=\cA^*_c(H)$ projects over the whole $M$ and therefore it is an invariant Lipschitz Lagrangian graph.
\end{lemma}

\begin{Proof}
Let $c\in \rH^1(M;\R)$ and consider a connected component $U^*$ of $\cM^*_c(H)$ and let $U=\pi(U^*)$, where $\pi: \rT^*M\longrightarrow M$ denotes the canonical projection along the fibers. Obviously $U$ is closed ($\pi|\cM^*_c(H)$ is a bi-Lipschitz homeomorphism, see Mather's graph theorem \cite{Mather91}). 
We would like to show that $U$ is also open; if this is the case, then $U=M$ (because of the connectedness of $M$) and therefore, using the graph property, $\cM^*_c(H)=\cA^*_c(H)$. The fact that it is a Lipschitz Lagrangian graph follows from weak KAM theory's characterization of the Aubry set in terms of subsolutions of Hamilton-Jacobi equation (see (\ref{Aubry}) and Remark \ref{Remark11}).

To show that $U$ is open it is sufficient to use that $U^*$ is contained in an $n$-dimensional graph and it is invariant under the action of $n$ independent vector fields that commute with $X_H$.
Denote by $F_1,\ldots, F_n$ the $n$ independent integrals of motion and by $\Phi_{F_1}, \ldots, \Phi_{F_n}$  their respective flows. These Hamiltonian flows are well-defined since we are assuming that the integrals of motion $F_i$'s are $C^2$ (or at least $C^{1,1}$).
Let us consider a point $(x_0,p_0)\in U^*$ and define $\G_0:=\{(x_0,p_0)\}$. Consider now the evolution of $\G_0$ under the flow $\Phi_{F_1}$. Let $\e_1>0$ be sufficiently small and define $\G_1:=\{\Phi^t_{F_1}(x_0,p_0),\; |t|\leq \e_1\}$. In the same way, taking $\e_2>0$ sufficiently small, one can define 
$\G_2:=\{\Phi^t_{F_2}(x,p),\;(x,p)\in \G_1 \;\mbox{and}\; |t|\leq \e_2\}$. This set is fibered by $\G^{(x,p)}_2:=\{\Phi^t_{F_2}(x,p),\; |t|\leq \e_2\}$ for each $(x,p)\in \G_1$. Using the independence of the vector fields $X_{F_1}$ and $X_{F_2}$ and the smallness assumptions on $\e_1$ and $\e_2$, one can assume that these fibers are disjoint and intersect $\G_1$ only once. In the very same way, one can define $\G_3,\ldots,\, \G_k, \, \ldots ,\G_n$ and choose $\e_3,\ldots, \e_n$ sufficiently small so that each $\G_k$ is fibered over $\G_{k-1}$ and all fibers are disjoint and intersect $\G_{k-1}$ only once. Let us now consider $\e<\min \{\e_1,\ldots,\e_n\}$ and define
\beqano
\Theta : [-\e,\e]^n &\longrightarrow& U^* \\
(t_1,\ldots,t_n) &\longmapsto&  \Phi^{t_n}_{F_n}\big(\ldots \ldots \Phi^{t_2}_{F_2}\big(\Phi^{t_1}_{F_1}(x_0,p_0)\big) \ldots\big)\,.
\eeqano
$\Theta$ is clearly continuous (because of the continuity of the flows) and it is injective. 
This allows us to conclude that there exists an open neighborhood of $x_0$ in $U$ and consequently each $x_0$ is contained in the interior of $U$. 
This concludes the proof.

Alternatively, one can observe that the map $\Theta$ is $C^1$ and has injective derivative at $(0,\ldots,0)$, therefore its image (which is contained in the Mather set) contains a piece of an $n$-dimensional $C^1$ submanifold around $(x_0,p_0)$. Using Brouwer's theorem \cite{Brouwer}, one can conclude that the projection of the Mather set on $M$ is open, hence everything.
\end{Proof}

\begin{remark}\label{remarkregularity}
It is also possible to deduce from the above proof (see in particular the last part), that the Mather set is a $C^k$ manifold $N$, if all $F_i$'s are $C^k$, $k\geq 2$. In particular, it is a $C^{k-1}$ graph over the base. This follows from the fact that $N$ is a Lipschiz graph (Mather's graph theorem \cite{Mather91}) and therefore it is transversal to the fibers of the cotangent bundle. 
\end{remark}

Regularity results for the Mather (or Aubry) set can be also deduced using the following lemma, via the implicit function theorem.

\begin{lemma}\label{lemma123}
If $F: \rT^*M \rightarrow \R$ is a Hamiltonian and $\L$ is a Lipschitz Lagrangian graph invariant under the flow of $F$, then $F$ is constant on $\L$.
\end{lemma}

\begin{Proof}
Let $\L$ be the graph of a closed $1$-form $\eta: M \rightarrow \rT^*M$. Recall that 
$\eta$ is differentiable almost everywhere and that at differentiability points $\rT_{\eta(x)}\Lambda$ is a Lagrangian subspace.
Let $x_1, x_2  \in M$  and let $\g$ be an absolutely continuous curve in $M$ connecting them, such that $\eta$ is differentiable almost everywhere along $\g$ (with respect to the $1$-dimensional Lebesgue measure). In fact, one can consider the family $\cC_{x_1,x_2}$ of all absolutely continuous curves $\xi$ that connect $x_1$ to $x_2$ and whose length is, for instance, less than  $2\, {dist}(x_1,x_2)$. If one denotes by $U$ the set spanned by these curves, it is easy to see that  $U\setminus \{x_1,x_2\}$ is open and therefore $\eta$ is differentiable almost everywhere in $U$. Applying Fubini's theorem one can deduce the existence in $\cC_{x_1,x_2}$ of a curve $\g$ as above. \\
Let us now consider the lift of $\g$ onto $\Lambda$, \ie $\G:= \eta \circ \g$. 
This curve $\G$ is differentiable almost everywhere and at each differentiability point $t_0$ we have $\frac{d}{dt} F(\G(t_0))= dF(\G(t_0))\cdot \dot{\G}(t_0)$, with $\dot{\G}(t_0) \in \rT_{\G(t_0)}\L$. Using the fact that $\rT_{\G(t_0)}\L$ is a Lagrangian subspace and that $X_F$ is tangent to $\L$ at all points (since $\L$ is invariant), we can conclude that
$\frac{d}{dt} F(\G(t_0))= dF(\G(t_0))\cdot \dot{\G}(t_0) = \omega(X_F(\G(t_0)), \dot{\G}(t_0)) = 0$. Therefore, $\frac{d}{dt} F(\G(t))=0$ almost everywhere  and integrating along $\G$ we obtain that $F(\eta(x_1))=F(\eta(x_2))$.
\end{Proof}

We shall also need the following classical lemma (for its proof we refer the reader to \cite[Section 49]{Arnoldbook})

\begin{lemma}\label{step3}
Let $N^n$ be a compact connected differentiable $n$-dimensional manifold, on which we are given $n$ pairwise commuting $C^1$ vector fields, linearly independent at each point. Then $N^n$ is diffeomorphic to an $n$-dimensional torus. 
\end{lemma}

We can now prove our weak version of Liouville's theorem.\\

\noindent{\bf Proof} [{\bf Theorem \ref{teolio} (Weak Liouville's theorem)}].
i) The existence of these smooth Lagrangian graphs follows from Lemmata \ref{weakLiouville} and \ref{lemma123}. In fact, if we denote by 
$F_1,\ldots,F_n$ the $n$ independent integrals of motion, then
for each $c\in \rH^1(M;\R)$ we have that 
$\cM^*_c(H)=\cA^*_c(H) $ coincides with a connected component of $\{F_1=a_1,\ldots, F_n=a_n\}$ for some $(a_1,\ldots,a_n)\in \R^n$.
Therefore these Lagrangian graphs are smooth and each of them supports an action-minimizing measure of full support and with a certain rotation vector $h(c)\in \rH_1(M;\R)$ (see \cite{Mather91} for a precise definition of rotation vector). Observe, in fact,  that it is always possible to find a $c$-action minimizing measure, whose support is the whole Mather set $\widetilde{\cM}_c(L)$ (we consider the Lagrangian setting on $\rT M$ and then use the Legendre transform $\cL_L$ to push everything forward to $\rT^*M$). 
In fact, since the space of probability measures on ${\rm T}M$ is a separable metric space, one can take
a countable dense set $\{\m_n\}_{n=1}^{\infty}$ of $c$-action minimizing measures  and 
consider the new measure $\tilde{\m} = \sum_{n=1}^\infty \frac{1}{2^n}\m_n$. This measure is still invariant and $c$-action minimizing  and ${\rm supp}\,\tilde{\m} = \widetilde{\cM}_c(L)$. The rotation vector $h(c)$ will be the rotation vector of such a measure $\tilde{\m}$.
Moreover, because of the graph property, the Mather set corresponding to this rotation vector $h(c)$ must also coincide with $\L_c$, \ie ${\cM^{h(c)*}}(H) = \L_c$ (see \cite{Mather91} for a definition of this set).

{ii)} From the fact that each $\L_c$ coincides with ${\cM^{h(c)*}}(H)$ for some $h(c)\in \rH_1(M;\R)$, we can deduce that these Lagrangian graphs cannot intersect. In fact, if $\L_c \cap  \L_{c'} \neq 0$ then
${\cM^{h(c)*}}(H) \cap \cM^*_{c'}(H) \neq 0$, but this would imply that ${\cM^{h(c)*}}(H) \subseteq \cM^*_{c'}(H)$ and therefore $\L_c \subseteq \L_{c'}$ (see for instance  \cite{FGS}). Since they are both graphs over $M$, then 
 $\L_c = \L_{c'}$ and necessarily  $c=c'$ (they must have the same cohomology class).
Moreover, if $\L$ is another invariant (Lipschitz) Lagrangian graph with cohomology class $c$, then it must contain the Aubry set $\cA^*_c(H) = \L_c$ (see (\ref{Aubry}) and Remark \ref{Remark11}); therefore $\L$ must coincide with $\L_c$.

iii) Let us now consider the union of these Lagrangian graphs, namely $\L^*:=\cup_{c\in\rH^1(M;\R)} \L_c$. 
Clearly $\L^*$ is closed. In fact, $\Lambda_c$ are equi-Lipschitz graphs, as it can be deduced from Mather's graph theorem \cite[Theorem 2]{Mather91} (even if not explicitly stated, the Lipschitz constant $C$ can be chosen to be the same for cohomology classes in a given compact region - see \cite[Lemma on page 186]{Mather91}). Alternatively, one can say that these $\L_c$ are the graphs of the differentials of classical solutions of Hamilton-Jacobi, which are locally equi-$C^{1,1}$, hence  their differentials are locally equi-Lipschitz (see \cite{Fathibook}). \\
From Proposition \ref{step2}, it follows that the integrals of motion are in involution on $\L^*$. 
Therefore, using Lemma \ref{step3} (and Remark \ref{Rem456}), one can deduce that all 
$\L_c \subseteq \mbox{Int}(\L^*)$ are diffeomorphic to $n$-dimensional tori. Proceeding as in the usual proof of Liouville's theorem (see for instance \cite[Sections 49--50]{Arnoldbook}) one can show that  the motion on each $\L_c$ is conjugate to a rotation on $\T^n$ (with rotation vector $h(c)$), \ie in an open neighborhood of $\L_c$ there exists a symplectic change of coordinates that transforms $\L_c$ into $\T^n \times \{0\}$ and the motion on it into a rigid rotation on $\T^n\times \{0\}$.\\
If $\L^*$ is open then it coincides with $\rT^*M$ (because of the connectedness of $\rT^*M$) and hence the system would be integrable in the classical sense; in fact, in this case the integrals of motion would be in involution everywhere.\\
\qed
\vspace{10 pt}

From the properties of these Lagrangian graphs, in particular the fact that they are disjoint, we can also conclude the following regularity result for $\a_H$ and its convex conjugate $\b_H$, as pointed out in Remark \ref{Remdopoteolio}.

\begin{corollary}\label{Coralphabeta}
If $H$ is a weakly  integrable Tonelli Hamiltonian, then Mather's $\a$-function $\a_H$ is strictly convex and $\b_H$ is $C^1$. 
\end{corollary}

\begin{Proof}
It is easy to deduce that $\a_{H}$ is strictly convex from the disjointness of the above Lagrangian graphs (\ie the Aubry sets). In fact, suppose by contradiction that there exist $\l\in (0,1)$ and $c,c'\in \rH^1(M;\R)$ such that
$\a_H(\l c + (1-\l)c') = \l \a_H(c) + (1-\l)\a_H(c')$ and let $\m_{\l}$ be a $(\l c + (1-\l)c')$ - action minimizing measure. Let us denote by $\eta_c$ and $\eta_{c'}$ two closed $1$-forms with cohomology classes, respectively, $c$ and $c'$. Then:
\beqano
-\a_H(\l c + (1-\l)c') &=& \int \big(L - \l \eta_c - (1-\l)\eta_{c'} \big) d\m_{\l} = \\
&=& \l \int \big(L - \eta_c\big) d\m_{\l}  + (1-\l) \int \big(L - \eta_{c'} \big) d\m_{\l}  \geq \\
&\geq& - \l \a_H(c) - (1-\l)\a_H(c') = -\a_H(\l c + (1-\l)c'),
\eeqano
therefore all above inequalities are equalities and this implies that $\m_{\l}$ is also 
$c$-action minimizing and $c'$-action minimizing. Obviously this contradicts the disjointness of the Mather and Aubry sets (\ie the $\L_c$'s).
As far as the differentiability of $\b_H$ is concerned, it also follows from the disjointness of the $\L_c$'s: if $c,c' \in \dpr \b_H(h)$ for some $h\in \rH_1(M;\R)$, then the corresponding Mather sets $\cM^*_c$ and $\cM^*_{c'}$ would contain ${\cM^{h}}^*$ and therefore $\L_c$ and $\L_{c'}$ would intersect.\\
Observe that the differentiability of $\beta_H$ can be also deduced as a consequence of the strict convexity of $\a_H$ and of basic properties of the Fenchel-Legendre transform.
\end{Proof}

Now, let us prove Theorem \ref{teolio2}.\\

\noindent{\bf Proof} [{\bf Theorem \ref{teolio2}}]. Let us identify $\rH^1(M;\R)$ with $\R^n$. For each cohomology class $c\in \R^n$ let us consider the unique Lagrangian graph $\L_c = \{(x, c+ du_c):\; x\in M\}$ given by Theorem \ref{teolio}, where $u_c: M \rightarrow \R$ is a smooth function. In particular, it follows from Lemma \ref{lemma123}  that $\L_c \subseteq \{F_1=a_1(c),\ldots, F_n=a_n(c)\}$ for some $\vec{a}(c)=(a_1(c),\ldots,a_n(c))\in \R^n$. This allows us to define the following function:
\beqano
\cF : \R^n &\longrightarrow& \R^n\\
 c &\longmapsto& \vec{a}(c)=(a_1(c),\ldots, a_n(c)).
\eeqano
This function is clearly well-defined (because of the uniqueness of the $\L_c$'s). We want to show that it is also continuous, actually it is locally Lipschitz. Let $K$ be a compact subset of $\R^n$. It is easy to check that 
$\cup_{c\in K} \L_c$ is contained in a compact region of $\rT^*M$ and let us denote by 
$C(K)>0$ a common Lipschitz constant for $F_1,\ldots, F_n$ in such a region . Let us now consider $c,c'\in K$ and observe that there exists at least one point $x_0 \in M$ such that $d_{x_0}u_c=d_{x_0}u_{c'}$; in fact, the function $u_c-u_{c'}$ is a smooth function on a compact manifold, hence it has critical points. Then:
\beqano
\|\cF(c) - \cF(c')\|_{\infty} &=& \|\vec{a}(c)-\vec{a}(c')\|_{\infty} \leq\\
&\leq& \max_{i=1,\ldots, n} \|
F_i(x_0,c+d_{x_0}u_c) - F_i(x_0,c'+d_{x_0}u_{c'})\| \leq  \\
&\leq& C(K) \|(c+d_{x_0}u_c) - (c'+d_{x_0}u_{c'})\| =\\
&=& C(K) \|c-c'\|\,.
\eeqano

Let us show that $\cF$ is locally injective. Since $\L_c$ is compact and $F=(F_1,\ldots, F_n)$ is a submersion constant on $\L_c$, we can find $\cN\rightarrow \L_c$ a tubular neighborhood of $\L_c$ and $\e>0$, such that on each fiber of $\cN$, the map $F$ is a diffeomorphism onto the ball $B(\vec{a}(c),\e)$, where $\vec{a}(c)\in\R^n$ is the value of $F$ on $\L_c$.
This implies that the levels of $F$ intersect $\cN$ in a connected set diffeomorphic to $\L_c$. For $c'$ close to $c$, we have $\L_{c'}\subset \cN$ and hence $\L_{c'}$ is equal to the intersection of
$F^{-1}(\cF(c'))$ with $\cN$.
Observe that so far we have not used that $\dim H^1(M;\R)=n$: the above argument makes perfectly sense even if $\dim H^1(M;\R)<n$.

Now, if we assume that   $\dim H^1(M;\R)=n$, then local injectivity of $\cF$ implies that it is a local homeomorphism. Then, the argument above shows that a neighborhood of $\L_c$ is contained in $\L^*:=\cup_{c}\L_c$. From which integrability follows.
 
 \qed

\section{Proof of Lemma \ref{maintheorem} }\label{secproof}

The proof of Lemma \ref{maintheorem}  will follow from the following lemmata. First of all, let us recall this classical result in Hamiltonian dynamics, whose proof can be found, for instance, in \cite{Arnoldbook}.\\

\begin{lemma}\label{lemma1}
Let $H$ and $F$ be two Hamiltonians on $\rT^*M$. Then:
\beqano
H \;\mbox{is constant on the orbits of}\; F\;  &\Longleftrightarrow& \;F \;\mbox{is constant on the orbits of}\; H\\
&\Longleftrightarrow& \{H,F\}=0\,.
\eeqano
Moreover, if $\{H,F\}=0$ then the two flows $\Phi_H$ and $\Phi_F$ commute, \ie
$\Phi^t_H \circ \Phi^s_F = \Phi^s_F \circ \Phi^t_H$ for all $s,t\in \R$.\\
\end{lemma}

\begin{remark}\label{Rem456}
Observe that if $\{H,F\}=0$ in an open region, then $\Phi_H$ and $\Phi_F$ will commute in that region. This condition is sufficient, but not necessary. One can show that $[X_H,X_F] = - X_{\{H,F\}}$, 
where $[\cdot,\cdot]$ denotes the commutator between two vector fields and $X_G$ the Hamiltonian vector field associated to a Hamiltonian $G$. Therefore it is easy to  check that
$\Phi_H$ and $\Phi_F$ commute if and only if $\{H,F\}$ is locally constant (see \cite{Arnoldbook})\\
\end{remark}

\begin{lemma}\label{lemma2}
If $H$ and $F$ are two commuting Hamiltonians and $\m$ is an invariant measure of $\Phi_H$, then, for each $t\in \R$, ${\Phi^t_F}_*\m$ is still  $\Phi_H$-invariant.
\end{lemma}

\begin{Proof}
Let $\tilde{\m}_t:={\Phi^t_F}_*\m$ denote the push-forward of $\m$. We want to show that for each $s\in \R $, 
${\Phi^s_H}_*\tilde{\m}_t = \tilde{\m}_t$. In fact, since ${\Phi^t_H}$ and ${\Phi^s_F}$ commute (because $\{H,F\}=0$) and $\m$ is $\Phi_H$-invariant (\ie ${\Phi^s_H}_*\m =\m$), we have:
\beqano
{\Phi^s_H}_*\tilde{\m}_t &=&  {\Phi^s_H}_*\left({\Phi^t_F}_*{\m}\right) = {\Phi^t_F}_*\left({\Phi^s_H}_*{\m}\right) = {\Phi^t_F}_* {\m} = \tilde{\m}_t\,.
\eeqano
\end{Proof}

In the following, it will be convenient to consider this characterization of 
$c$-action minizimizing measures, that was proved in \cite{FGS}. First of all, we need to recall the definition of $c$-subcritical Lagrangian graph (also introduced in \cite{FGS}).

\begin{definition}[\bf $c$-subcritical Lagrangian graph]\label{Def2} Given a Lipschitz Lagrangian graph $\Lambda$ with cohomology class $c$, we shall say that
$\Lambda$ is $c$-{\it subcritical} 
for a Tonelli Hamiltonian $H$, if 
$$\Lambda \subset \{(x,p) \in {\rm T}^*M:\; H(x,p)\leq \a_H(c)\},$$
where $\a_H :{\rm H}^1(M;\R) 
\longrightarrow \R $ is Mather's $\a$-function associated to $H$. 
We shall call its {\it critical} part: $\Lambda_{crit}=\{(x,p)\in \Lambda:\; H(x,p) = \a_H(c)\}$.\\
\end{definition}

\begin{remark}\label{Remark1}
The interest in such graphs comes from the fact that 
$\{H(x,p)\leq \a_H(c)\}$
is the smallest energy sublevel of $H$ containing Lipschitz Lagrangian graphs of cohomology class $c$. It follows from the results in 
\cite{Fathibook, bernardc11}, that they do always exist. 
As we have already recalled in (\ref{Aubry}),
$\cA^*_c(H)$ can be characterized as the intersection of all the $c$-subcritical Lagrangian graphs of $H$ or equivalently of all their critical parts. Therefore, the critical part of these Lagrangian graphs 
is always non-empty and contains the Aubry set $\cA^*_c(H)$. 
Moreover, there always exists a $c$-subcritical Lagrangian graph $\tilde{\Lambda}$ such that $\tilde{\Lambda}_{crit} = \cA^*_c(H)$ (see again \cite{Fathibook}).
\end{remark}

In \cite{FGS} we proved the following characterization of $c$-action minimizing measures.

\begin{lemma} \label{lemmadellacontesa}
Let $\m$ be an invariant probability measure for a Tonelli Hamiltonian $H$ on ${\rm T}^*M$.
$\m$ is a $c$-action minimizing measure if and only if ${\rm supp}\,\m$ is contained in the critical part of a $c$-subcritical Lagrangian graph of $H$. In particular, 
any invariant probability measure $\m$, whose support
is contained in an invariant Lagrangian graph of $H$ with Liouville class $c$,
is $c$-action minimizing.
\end{lemma}

We have now recalled all the needed ingredients for the proof of Lemma \ref{maintheorem}.\\

\noindent{\bf Proof [Lemma \ref{maintheorem}].} 
First of all, let us observe that it is enough to show the results for $|t|\leq \e_0$, for some sufficiently small $\e_0>0$.\\
(\rm i) From Lemma \ref{lemma2}, we know that $\tilde{\m}_t:={\Phi^t_F}_*\m$ is still invariant under the action of $\Phi_H$. We need to show that it is still $c$-action minimizing. 
In the light of Lemma \ref{lemmadellacontesa}, it will be sufficient to prove that $\supp \tilde{\m}_t$ is contained in the critical part of a $c$-subcritical Lagrangian graph of $H$. Hence, let $\Lambda$ be any $c$-subcritical Lagrangian graph of $H$. This graph contains $\supp \m$ in its critical part because $\m$ is $c$-action minimizing and $\cM^*_c(H) \subseteq \cA^*_c(H)$ (see Remark \ref{Remark1}). Consider now $\tilde{\Lambda}:=\Phi^{t}_F(\Lambda)$; this is also a Lagrangian manifold of cohomology class $c$ (since $\Phi^{t}_F$ is an exact symplectomorphism) and using that $\Lambda$ is Lipschitz it is easy to see that there exists $\e_0=\e_0(F,\Lambda)>0$ such that
$\tilde{\Lambda}$ is still a graph for all $|t|\leq\e_0$. Moreover, it follows from Lemma \ref{lemma1} that $\tilde{\Lambda}$ is still $c$-subcritical for $H$.
In order to conclude the proof, it is enough to observe that  $\supp \tilde{\m}$ is contained in the critical part of $\tilde{\Lambda}$ (it follows also from Lemma \ref{lemma1}).\\
{\rm (ii)} The invariance of the Mather set follows immediately from (i) and its definition (it is the  union of all $c$-action minimizing measures).
To prove that the $\cA_c^*(H)$ is invariant under $\Phi^t_F$, we proceed exactly as before, 
using the fact that, as recalled in Remark \ref{Remark1}, $\cA^*_c(H)$ can be obtained by intersecting all $c$-subcritical Lagrangian  graphs of $H$. In particular, there exists a $c$-subcritical Lagrangian graph $\Lambda$ such that $\cA^*_c(H)=\Lambda_{crit}$.
For $|t|\leq \e_0$, as in (i), let us consider $\tilde{\Lambda}:=\Phi^{-t}_F(\Lambda)$. As we have already observed, this is still a $c$-subcritical Lagrangian graph of $H$ and therefore $\cA^*_c (H) \subseteq \tilde{\Lambda}_{crit}$; moreover, it follows from Lemma \ref{lemma1} that $\tilde{\Lambda}_{crit}=\Phi^{-t}_F (\Lambda_{crit})$. This is enough to conclude the proof. In fact:
\beqano
\Phi^t_F\left(\cA^*_c(H)\right) \subseteq \Phi^t_F\left(\tilde{\Lambda}_{crit}\right) = \Lambda_{crit}= \cA^*_c(H)\,.
\eeqano
Furthermore, it is clear from the proof that each connected component of these sets is mapped into itself.
\qed
\vspace{15 pt}

To conclude let us observe that part $(i)$ of Lemma \ref{maintheorem} can be also showed in a more direct way, without passing through Lemma \ref{lemmadellacontesa}. However, in order to show the symplectic invariance of the Aubry set, extra tools are necessary.  \\

\begin{lemma}
Let $\m$ be an invariant probability measure for a Tonelli Hamiltonian $H$ on $\rT^*M$ and $\Phi: \rT^*M \longrightarrow \rT^*M$ an exact symplectomorphism that preserves $H$, \ie $H\circ \Phi = H$. Then:  
$$
\int \left[p\frac{\dpr H}{\dpr p}(x,p) - H(x,p)\right] d \m =
\int \left[p\frac{\dpr H}{\dpr p}(x,p) - H(x,p)\right] d \Phi_*\m.
$$
\end{lemma}

Recall that a symplectomorphism $\Phi: \rT^*M \longrightarrow \rT^*M$ is said to be exact if $\Phi_*(p dx) - pdx$ is an exact $1$-form.

\begin{remark}
Part $(i)$ of Lemma \ref{lemmadellacontesa} simply follows choosing $\Phi=\Phi_F^{t}$ and using the definition of action-miniziming measures (see section \ref{secmain}) and the relation between Hamiltonian and Lagrangian (Fenchel-Legendre transform).
\end{remark}

\begin{Proof}
If we denote by $\l(x,p)$ the Liouville form $pdx$ and by $X_H(x,p)$ the Hamiltonian vector field, then:
\beqano
\int \left[p\frac{\dpr H}{\dpr p}(x,p) - H(x,p)\right] d \m = 
\int \Big(\l(x,p)[X_H(x,p)] - H(x,p)\Big) d \m\,.
\eeqano
Therefore, using that $|\det D \Phi| =1$  
and $\Phi^*\l - \l = df$ (since $\Phi$ is an exact symplectmorphism) and that $H$ is preserved by $\Phi$ we obtain: 
\beqano
&& \int \left[p\frac{\dpr H}{\dpr p}(x,p) - H(x,p)\right]\;d\Phi_*\m = 
\int \Big(\l(x,p)[X_H(x,p)] - H(x,p)\Big)\;d \Phi_*\m = \\
&& \; =\; \int \Big(\Phi^*\l(x',p')[X_{H\circ \Phi}(x',p')] - H(\Phi(x',p'))\Big)\;d \m = \\
&& \; =\; \int \Big(\l(x',p') [X_H(x',p')] - H(x',p')\Big)\;d \m +  \int df(x',p')[X_H(x',p')] \,d\m=\\
&& \;= \int \Big(p'\frac{\dpr H}{\dpr p}(x',p') - H(x',p')\Big)\;d \m. 
\eeqano
In the last equality we used that 
$\int df(x',p')[X_H(x',p')] \,d\m =0$, as it follows easily from the invariance of $\m$. In fact, let us assume that $\m$ is ergodic (otherwise consider each ergodic component); using the ergodic theorem and the compactness of $M$, we obtain that for a generic point $(x_0,y_0)$ in the support of $\m$:
\beqano
\int df(x',p')[X_H(x',p')] \,d\m &=& \lim_{N\rightarrow +\infty} \frac{1}{N} \int_0^N
df(\Phi_H^t(x_0,p_0))[X_H(\Phi_H^t(x_0,p_0))]\,d t = \\
&=& \lim_{N\rightarrow +\infty} \frac{f(\Phi_H^N(x_0,p_0)) - f(x_0,p_0)}{N} =0\,. 
\eeqano

Alternatively, without using the ergodic theorem, one  can observe that
$$
f(\Phi_H^t(x',p')) - f(x',p') = \int_0^t df(\Phi_H^t(x',p'))(X_H(\Phi_H^t(x',p')))\,dt
$$
and integrate with respect to $\mu$. Then, the result follows applying Fubini's theorem and using the invariance of $\m$ by $\Phi^t_H$.
\end{Proof}


\appendix

\section{Examples of weakly integrable Hamiltonians} \label{App1}

In this appendix we would like to exhibit an interesting class of weakly integrable Tonelli Hamiltonians, namely the geodesics flows associated to left-invariant Riemannian metrics on Lie groups. These flows have been extensively studied in classical mechanics; 
 in fact they represent, in some sense, a natural generalization of Eulerian motions of a rigid body   this point of view can be dated back at least to  Henri Poincar\'e's article \cite{Poincare}).\\
In the following, we shall mainly follow \cite[Appendix 2]{Arnoldbook}, but we also refer the reader to \cite[Chapter VI, 1.B]{Arnoldbook2} and \cite{Fedorov, Poincare} for more details on the subject.

Let us consider a compact Lie group $G$ and denote by $\fg$ its Lie algebra, \ie the tangent space to the group at the identity element. A Riemannian metric on $G$ is said to be {\it left-invariant} if it is preserved by all left translations $L_{g}: h \mapsto gh$, \ie the derivatives of all left translations map every vector in $\rT_hG$ into a vector of the same length in $\rT_{gh}G$. Obviously it is sufficient to specify the metric at one point of the group, for instance the identity element $e\in G$, and therefore there are as many left-invariant metrics on $G$ as there are {\it euclidean structures} on $\fg$, \ie  symmetric positive definite operators from the algebra to its dual space: $A: \fg \rightarrow \fg^*$.\\
We would like to consider the {\it geodesic flow} associated to such a left-invariant metric space $(G,A)$. Sometimes this flow is also referred to as {\it motion of a generalized rigid body with configuration space} $G$. In classical mechanics, in fact,  $G=SO(3)$ - \ie the group of rotation of a $3$-dimensional euclidean space -  can be regarded as the configuration space of a rigid body fixed at a point (while the Lie algebra $\fg$ represents the $3$-dimensional space of angular velocities of all possible rotations) and the motion of the body can be described by curves $g=g(t)$ in $G$, that correspond to geodesics of a left-invariant metric. 
Let us see what happens with a generic Lie group $G$. First of all observe that the euclidean structure $A$ can be extended to all fibers:
\beqano
{A_{g}}: \rT_g G &\longrightarrow& \rT_g^*G\\
\dot{g} &\longmapsto& L^*_{g^{-1}}A{L_{g^{-1}}}_*\dot{g},
\eeqano
where $L^*_{g^{-1}}$ and ${L_{g^{-1}}}_*$ are respectively the maps induced by the left translation $L_{g^{-1}}$ to the cotangent and the tangent space of $G$. The operator $\tilde{A}(g,\dot{g}):=(g,{A}_{g}(\dot{g}))$ is called the {\it moment of intertia operator}. This allows us to define the Lagrangian  associated to the geodesic flow on $(G,A)$:
\beqano
L : \rT G &\longrightarrow \R&\\
(g,\dot{g}) &\longmapsto& \frac{1}{2} \langle {A}_g\dot{g}, \dot{g} \rangle\,;
\eeqano
in other words, this represents the {\it kinetic energy} of the system. It is easy to check that the associated Hamiltonian is:
\beqano
H : \rT^* G &\longrightarrow \R&\\
(g, p) &\longmapsto& \frac{1}{2} \langle p, {A}^{-1}_g p \rangle\,.
\eeqano
Such a Hamiltonian is a Tonelli Hamiltonian on $\rT^* G$. 
Let us study its integrals of motion.
First, observe that each angular momentum $p\in \rT_g^*G$ can be  carried to $\fg^*$ by both left and right translation: $p_b=L^*_gp$ and $p_s=R^*_gp$ ($R_{g}$ denotes the right translation $h \mapsto hg$); in classical mechanics these two vectors are called respectively {\it angular momentum relative to the body} and {\it relative to the space}. Euler showed (in the case $G=SO(3)$, but the same proof works for a general $G$) that the motions of these two angular momenta satisfy the following equations (known as {\it Euler's equations for the rigid body}):
$$
\frac{d p_s}{dt}=0 \qquad {\rm and }\qquad \frac{d p_b}{dt}= ad^*_{A^{-1}p_b}p_b\,,
$$
where $ad^*_{\x}: \fg^* \rightarrow \fg^*$ denotes the so-called {\it co-adjoint representation of the group}. In particular, the second equation determines a flow $\varphi_t: \fg^*\rightarrow \fg^*$, which describes the motion of the angular momentum relative to the body (observe that it does not depend on the position of the body in the space). If one defines the map $\pi: \rT^*G \rightarrow \fg^*$ given by $\pi(g,p)=L^*_gp$, it is not difficult to see that this map is a {\it factorization} of the Hamiltonian flow $(\rT^*G,\Phi_H^t)$ over the flow $(\fg^*,\varphi^t)$, \ie the following diagram commutes:
$$\xymatrix{
{\rT^* G} \ar@{->}[d]_{\pi} \ar@{->}[r]^{\Phi^t_H}  & {\rT^*G}\ar@{->}[d]^{\pi} 
\\  
{\fg^*} \ar@{->}[r]_{\varphi^t}  & {\fg^*} \\}
 $$
From this and the conservation law for the vector of angular momentum relative to the space (in particular each of its components is conserved), we obtain a set of integrals of motion for $H$. Observe that to each element of the Lie algebra $\fg$ corresponds a linear functional on the space $\fg^*$ and therefore an integral of motion of $H$. 
Of course, of all these integrals of motion at most $n$ can be functionally independent. For instance, one can take the ones obtained by $n$ linear functionals on $\fg^*$ which form a basis in $\fg$. This proves the following.

\begin{proposition}\label{propgruppilie}
All Tonelli Hamiltonians corresponding to geodesic flows associated to left-invariant Riemannian metrics on compact Lie groups are weakly integrable.
 \end{proposition}

\begin{remark}
Observe that in general these geodesic flows are not necessarily integrable in the sense of Liouville (in the form stated in Section \ref{intro}). The problem of the {\it non-integrability} in the sense of Liouville is extremely subtle and tricky, and may depend on the regularity class in which we are looking for the integrals of motion (analytic, smooth, {\it etc}...) or on which sense of independence we ask (independent everywhere, on an open dense set, {\it etc} ...). See for instance \cite{Butler, Butler2} and references therein.  
\end{remark}

\def\cprime{$'$}

\end{document}